\documentclass{ifacconf}

\usepackage{graphicx}      
\usepackage{natbib}        
\usepackage{amsmath}
\usepackage{amssymb}
\usepackage{mathrsfs}  
\usepackage{lipsum}
\usepackage{mathtools}
\usepackage{cuted}
\usepackage{caption}
\usepackage{subcaption}
\usepackage{algorithm}
\usepackage[noend]{algpseudocode}

\usepackage{color}

\newcommand{\ie}{\emph{i.e.}}
\newcommand{\eg}{\emph{e.g.}}

\newcommand{\E}{\mathbb{E}}

\newcommand{\abs}[1]{\left\vert #1 \right\vert}

\renewcommand{\t}{\theta}

\newcommand{\amin}[1]{\underset{#1}{\operatorname{arg\,min}}} 
\renewcommand{\sup}[1]{\underset{#1}{\operatorname{sup}}}

\renewcommand{\S}[3]{\sum_{\substack{#1 \\ #2}}^{#3}}

\renewcommand{\amin}[1]{\underset{#1}{\operatorname{arg\,min}}}

\renewcommand{\S}[3]{\sum_{\substack{#1 \\ #2}}^{#3}}

\newcommand{\pb}{\delta_{lpb,\varepsilon}}

\usepackage{mathtools}  
\mathtoolsset{showonlyrefs}

\begin{document}
\begin{frontmatter}

\title{A Unified Approach to Differentially Private Bayes Point Estimation \thanksref{footnoteinfo}} 

\thanks[footnoteinfo]{This work has been supported by the Swedish Research Council under contract number 2016-06079 (NewLEADS) and by the Digital Futures project EXTREMUM. The authors are with the Division of Decision and Control Systems, KTH Royal Institute of Technology, 100 44 Stockholm, Sweden (e-mails: blak@kth.se; crro@kth.se).}

\author[First]{Braghadeesh Lakshminarayanan} 
\author[First]{Cristian R. Rojas} 

\address[First]{Decision and Control Systems, KTH Royal Institute of Technology, 10044 Stockholm, Sweden (e-mails: blak@kth.se, crro@kth.se).}

\begin{abstract}                
Parameter estimation in statistics and system identification relies on data that may contain sensitive information. To protect this sensitive information, the notion of \emph{differential privacy} (DP) has been proposed, which enforces confidentiality by introducing randomization in the estimates. Standard algorithms for differentially private estimation are based on adding an appropriate amount of noise to the output of a traditional point estimation method. This leads to an accuracy-privacy trade off, as adding more noise reduces the accuracy while increasing privacy. In this paper, we propose a new Unified Bayes Private Point (UBaPP) approach to Bayes point estimation of the unknown parameters of a data generating mechanism under a DP constraint, that achieves a better accuracy-privacy trade off than traditional approaches. We verify the performance of our approach on a simple numerical example. 
\end{abstract}


\end{frontmatter}

\section{Introduction}
Parameter estimation deals with the problem of approximating the unknown parameters of a mathematical model that describes a given real phenomenon, using data collected from that phenomenon. This problem has been intensely studied in areas such as statistics~\citep{casella2021statistical}, system identification~\citep{Soderstrom-Stoica-89, Ljung:99}, and machine learning~\citep{Shalev-Shwartz-Ben-David-12}. 

An important subfield of parameter estimation is point estimation, where the goal is to approximate the unknown quantity by a single value.
Some of the most commonly used point estimators are Maximum Likelihood~\citep{casella2021statistical,LehmCase98}, the Method of Moments~\citep{Gourieroux-89}, the Minimum Mean Square Error estimator~\citep{van2004detection}, the Minimum Variance Unbiased estimator~\citep{Kay97}, and the Best Linear Unbiased estimator~\citep{BLUE}.

Given a point estimate, or in general a function of the data, it is possible to retrieve some information about the individual samples used to compute it.
For example, it was demonstrated in genomic studies~\citep{Homer} that, under some conditions, it is possible to identify whether or not the DNA sample of an individual was present in a dataset based on aggregate statistics. For this reason, the National Institute of Health (NIH) removed access to some aggregate statistics such as p-values and chi-squared statistics, which were once openly available~\citep{NIH}.
This and related concerns on confidentiality related to data handling has led to the development of point estimators subject to privacy considerations.




An important notion of privacy that is considered in the literature is \emph{differential privacy} (DP)~\citep{Dwork}. DP is a privacy constraint that can be imposed on algorithms in order to protect the sensitive information contained in their output. DP ensures that by seeing the output of an algorithm, almost no probabilistic inference can be made about the observations used by such algorithm to produce this output; the level of desired privacy can be tuned via a parameter $\varepsilon > 0$. A well known approach, known as the Laplace mechanism, enforces the DP constraint by adding a suitable amount of Laplace distributed random noise to the output of the algorithm.

Within automatic control, privacy has been well studied. For instance, \citep{sankar2011theory,Varodyan,9691832} consider information theory approaches to satisfy privacy constraints. DP has also been considered, \eg, in \cite{Pappas,wang2018differentially}.

One of the main issues that arise while enforcing DP in point estimation is the so-called accuracy-privacy trade off~\citep{wang2017differential,cao2020differentially}: Enforcing a higher level of privacy (for example, by adding Laplace noise of larger variance) reduces the accuracy of the estimator.
Most of the works that consider DP in parameter estimation rely on the Laplace mechanism to enforce DP. However, it is not clear if this mechanism achieves an optimal accuracy-privacy trade off.
%
In this paper, we provide an alternative approach that achieve an optimal accuracy-privacy trade off for Bayes point estimation by posing the problem of maximizing accuracy subject to a DP constraint as a convex optimization program.
In particular, our contributions are the following:
\begin{itemize}
    \item We formulate the problem of Bayes point estimation subject to a DP constraint as a convex optimization program;
    \item we provide an approach (UBaPP) to solve the above optimization program for the case where the parameter space and observations are discrete; 
    \item we demonstrate the advantage of our approach via a simple numerical example based on Bernoulli samples. 
\end{itemize}

The paper is organized as follows: Section~\ref{sec: Preliminaries} defines the notion of DP and reviews Bayes point estimation. In Section~\ref{sec: Setup}, we state the problem formulation, while in Section~\ref{sec: Approach} we propose our new approach (UBaPP). Then, we demonstrate our approach through a numerical example in Section~\ref{sec: Simulation}, and in Section~\ref{sec: Conclusion} we conclude the paper and discuss future work.

\section{Preliminaries} \label{sec: Preliminaries}

In this section, we formalize the problem of Bayes point estimation, and introduce DP. Then, we define the \emph{Laplace mechanism}, a procedure that enforces DP.  

\subsection{Bayes Point Estimation} \label{subsec: Bayesian Point Estimation}
Consider observing a physical process that generates independent and identically distributed (\emph{i.i.d.}) samples $x_1,\ldots,x_n \in \mathbb{R}^p$ at discrete time instants $1,\ldots,n$ respectively, according to a probability distribution $\mathbb{P}(\cdot \vert\theta)$ that is parameterized by an unknown parameter $\theta \in \Theta \subseteq \mathbb{R}^{d_\theta}$.
Let $\mathbf{x} := (x_1,\ldots,x_n)^T \in \mathbb{R}^{n \times p} =: \mathcal{X}$.
Let $S$ be a \emph{sufficient statistic}~\citep{LehmCase98} $S\colon \mathcal{X} \to \mathcal{Y} := \mathbb{R}^m$ where $\mathbf{y}=(y_1,\ldots,y_m)^T=S(\mathbf{x})$, with $m \leq n$. We call $\mathbf{y}$ \emph{observations} and $\mathcal{Y}$ the \emph{observation space}, and $\mathbf{x}$ \emph{input} and $\mathcal{X}$ the \emph{input space}.     

Given $\mathbf{x} \in \mathcal{X}$, the goal in point estimation is to construct an \emph{estimator} or \emph{decision rule}\footnote{In this section we focus on deterministic estimators, while in Section~\ref{sec: Approach} we extend these definitions to \emph{randomized} estimators.}, which is a mapping $\delta\colon \mathcal{Y} \rightarrow \Theta$, such that the \emph{risk}
\begin{equation} \label{eq: Bayes risk}
    R(\theta,\delta) := \mathbb{E} [L(\theta,\delta(S(\mathbf{x}))]
\end{equation}
is as small as possible, where $L\colon \Theta \times \Theta \rightarrow \mathbb{R}_0^{+}$ is a \emph{loss function}.

The expectation in \eqref{eq: Bayes risk} is taken with respect to the probability distribution of $\mathbf{x}$, \ie, $\mathbb{P}^n(\cdot\vert \theta)$, since the samples $x_1,\ldots,x_n$ are i.i.d.

The quantity $L(\theta,\delta(S(\mathbf{x})))$ measures the cost incurred in estimating the unknown parameter as $\delta(S(\mathbf{x}))$, whereas the true value of the parameter is $\theta$. Notice that $R(\theta,\delta)$ depends on $\theta$, which is unknown. Thus, in order to evaluate the performance of the estimator $\delta$, it is required to reduce $R(\theta,\delta)$ to a  function that depends only on $\delta$. For this purpose, Bayes point estimation~\citep{Kay97} assumes a \emph{prior} probability distribution $\pi(\cdot)$ over $\Theta$ and considers the average risk $\E_{\t \sim \pi}[R(\t, \delta)]$. Then, the notion of an ``optimal" Bayes decision rule is defined as
\begin{equation} \label{eq: Bayesian point estimate}
    \delta^*_\text{Bayes}(\mathbf{x}) := \amin{\delta \in \Delta_\mathbf{x}} \, \, \E_{\t \sim \pi(\cdot)}[R(\t, \delta)],
\end{equation}
where $\Delta_\mathbf{x} := \{\delta\colon \mathcal{Y} \to \Theta \}$.

If the loss function is the squared error
\begin{equation*}
    L(\theta,\delta(S(\mathbf{x}))) = (\theta - \delta(S(\mathbf{x})))^2,
\end{equation*}
then \eqref{eq: Bayesian point estimate} is the \emph{conditional mean estimate}~\citep{kailath2000linear}
\begin{equation}\label{eq: conditional_mean_estimate}
    \delta^*_\text{Bayes}(\mathbf{x}) = \mathbb{E}[\theta \vert S(\mathbf{x})],
\end{equation}
where the expectation is with respect to the posterior distribution of $\theta$ after observing $\mathbf{x}$.

We call \eqref{eq: Bayesian point estimate} a \emph{non-private Bayes point estimate}, since it does not consider any privacy constraints.
Similarly, we call \eqref{eq: conditional_mean_estimate} the \emph{non-private conditional mean estimate}.

\subsection{Differential Privacy} \label{subsec: Differential privacy}
To define DP, we first need to define the notion of \emph{neighbouring inputs}. For this, we need to introduce some notation.

Let $\mathbf{x} = (x_1,\ldots,x_n)^T$ and $\mathbf{x}' = (x'_1,\ldots,x'_n)^T$ be elements of $\mathcal{X} := \mathbb{R}^{n \times p}$, called the \emph{input space}.
Here, $x_i, \, x'_i \stackrel{i.i.d.}{\sim} \mathbb{P}$, for $i=1,\ldots,n$, where $\mathbb{P}$ is some probability distribution.

\begin{defn}\label{def: neighbouring inputs}
(Neighbouring inputs; \citealt{Dwork}). $\mathbf{x}$ and $\mathbf{x}'$ are called \emph{neighbouring inputs} if $d(\mathbf{x},\mathbf{x}') = 1$, where $d$ is the Hamming distance~\citep{Hamming}, \ie, if $x_i \neq x'_i$ for some unique $i \in \{1,\ldots,n\}$, and $x_j = x'_j$ for all $j \in \{1,\ldots,n\} \setminus \{i\}$.
\end{defn}

\begin{defn} \label{def: Differential privacy}
($\varepsilon$-Differential Privacy; \citealt{Dwork}). 
Let $\Theta$ be either a subset of $\mathbb{R}^{d_\theta}$ or $\mathbb{C}^{d_{\theta}}$, $d_{\theta} \in \mathbb{N}$. For $\varepsilon > 0$, a randomized algorithm $\mathscr{A}\colon \mathcal{X} \rightarrow \Theta$ is \emph{$\varepsilon$-differentially private} ($\varepsilon$-DP) if for each pair of neighboring inputs $\mathbf{x},\mathbf{x}' \in \mathcal{X}$ and $T \subseteq \Theta$ it holds that 
\begin{equation} \label{eq: Differential privacy}
\text{Pr}[\mathscr{A}(\mathbf{x}) \in T] \leq e^{\varepsilon} \, \text{Pr}[\mathscr{A}(\mathbf{x}') \in T]. 
\end{equation}
\end{defn}

\subsubsection{Interpretation} \hfill \\
Definition~\ref{def: Differential privacy} implies that for small values of $\varepsilon$, the probability distribution of the output $t = \mathscr{A}(\mathbf{x})$ of the algorithm is almost the same (up to a multiplicative constant $e^\varepsilon$) for two neighbouring inputs $\mathbf{x}$ and $\mathbf{x}'$. Therefore, by looking at the output $t$, it is difficult to infer whether $\mathbf{x}$ or $\mathbf{x}'$ is its corresponding input, since the distribution of $t$ is almost indistinguishable for $\mathbf{x}$ and $\mathbf{x}'$, thereby guaranteeing privacy. 

\subsection{Laplace Mechanism}
A standard approach to achieve $\varepsilon$-DP is the \emph{Laplace mechanism}~\citep{Dwork}, which we describe below. To this end, we need the notion of $l_1$-sensitivity. 

\begin{defn} \label{def: l1 sensitivity}
($l_1$-sensitivity; \citealt{Dwork}).
The \emph{$l_1$-sensitivity} of a function $g\colon \mathcal{X} \rightarrow \Theta$ is
\begin{equation} \label{eq: l1-sensitivity}
    \sigma_g := \sup{\mathbf{x},\mathbf{x}'\colon \, d(\mathbf{x},\mathbf{x}')=1} \, \, \vert \vert g(\mathbf{x}) - g(\mathbf{x}') \vert \vert_1, 
\end{equation}
where $\vert\vert \cdot \vert\vert_1$ denotes the $l_1$ norm~\citep{horn2012matrix}.
\end{defn}

\begin{defn} \label{def: Laplace Mechanism}
(Laplace Mechanism; \citealt{Dwork}).
Given a function $g\colon \mathcal{X} \rightarrow \Theta$, the \emph{Laplace mechanism} is a randomized algorithm that outputs the vector $J(\mathbf{x})$ whose $i^{th}$ component is distributed according to $[J(\mathbf{x})]_i \stackrel{i.i.d.}{\sim} \text{Lap}\left([g(\mathbf{x})]_i,\frac{\sigma_g}{\varepsilon}\right) $, for $i=1,\ldots,d_\theta$.
Here, $[g(\mathbf{x})]_i$ denotes the $i^{th}$ component of $g(\mathbf{x})$, and $\text{Lap}(a,b)$ is the Laplace distribution with probability density function (pdf)
\begin{equation}
f_{\text{Lap}}(z;a,b) = \frac{1}{2b} \exp\left({-\frac{\abs{x-a}}{b}}\right),
\end{equation}
where $b >0$.
\end{defn}

\medskip

\begin{rem}\label{rem: Laplace mechanism}
Notice from Definition~\ref{def: Laplace Mechanism} that $[J(\mathbf{x})]_i = [g(\mathbf{x})]_i + \eta_i$, where $\eta_i \stackrel{i.i.d.}{\sim} \text{Lap}\left(0, \frac{\sigma_g}{\varepsilon}\right)$, \ie, $J(\mathbf{x}) = g(\mathbf{x}) + \boldsymbol{\eta}$, where $\boldsymbol{\eta}=(\eta_1,\ldots,\eta_{d_\theta})^T.$ Hence, DP is enforced by explicitly randomizing the deterministic quantity $g(\mathbf{x})$ via the addition of a Laplace noise vector $\boldsymbol{\eta}$.
\end{rem}


\section{Problem Formulation} \label{sec: Setup}

\begin{figure}
    \centering
    \includegraphics[width=0.8\linewidth]{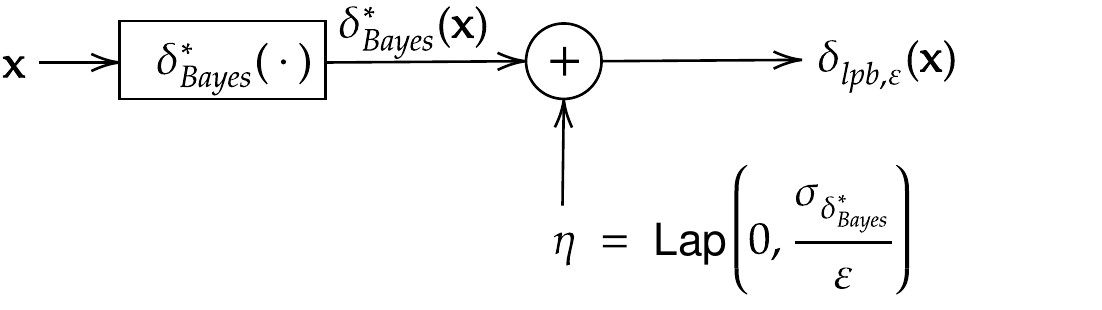}
    \caption{Cartoon illustrating the standard approach to obtain Bayes point estimate under DP.}
    \label{fig: naive private Bayes}
\end{figure}

In this paper we consider the problem of Bayes point estimation~\eqref{eq: Bayesian point estimate} under a DP constraint~\eqref{eq: Differential privacy}. In the standard Laplace mechanism approach, this is done in two steps. First, given input $\mathbf{x}$, a non-private Bayes point estimate is obtained using \eqref{eq: Bayesian point estimate}. Second,  Laplace noise is added to this non-private Bayes point estimate to obtain
\begin{equation}\label{eq: Laplace-private-Bayes}
    \pb(\mathbf{x}) = \delta_{\text{Bayes}}^*(\mathbf{x}) + \text{Lap}\left(0,\dfrac{\sigma_{\delta_{\text{Bayes}}^*}}{\varepsilon} \right).
\end{equation}
This quantity, called the \emph{Laplace Bayes Private Point} (LBaPP) estimate, satisfies the privacy constraint~\eqref{eq: Differential privacy}
(see Fig.~\ref{fig: naive private Bayes}.). Here,
\begin{equation}\label{eq: l1_sensitivity_Bayes}
    \sigma_{\delta_{\text{Bayes}}^*} = \sup{\mathbf{x},\mathbf{x}' \in \mathcal{X}\colon d(\mathbf{x},\mathbf{x}')=1} \, \, \vert \vert \delta_{\text{Bayes}}^*(\mathbf{x}) - \delta_{\text{Bayes}}^*(\mathbf{x}') \vert \vert_1.
\end{equation}
Although $\pb$ satisfies the privacy constraint, two important questions need to be addressed:
\begin{itemize}
    \item Is the Laplace noise addition ``optimal" in the sense of satisfying~\eqref{eq: Bayesian point estimate}?
    \item What if we do not have a closed form expression for $\delta_{\text{Bayes}}^*(\mathbf{x})$ that is required to compute $\sigma_{\delta_{\text{Bayes}}^*}$?
\end{itemize}
These two questions motivate the development, in the next section, of an alternative approach to find the optimal Bayes point estimator under a DP constraint.

\section{Proposed Approach}\label{sec: Approach}
Due to the definition of DP, we notice that
in order to impose privacy, we require the estimate to be randomized. We will allow the randomization of the estimator to be ``implicit" by replacing the \emph{deterministic} estimate $\delta(S(\mathbf{x}))$ with a \emph{randomized private Bayes estimate} $\delta_{p,\varepsilon}(\cdot\vert S(\mathbf{x}))$, which is a probability density function over the parameter space $\Theta$ for each input $\mathbf{x} \in \mathcal{X}$. This is formalized below.

Let $(\Theta, \mathcal{B}_\Theta)$ be a measurable space of parameters, $(\mathcal{X},\mathcal{B}_\mathcal{X})$ a measurable space of inputs endowed with a metric (e.g., the Hamming distance) $d$, $(\mathcal{Y},\mathcal{B}_\mathcal{Y})$ a measurable space of observations, $(q_\theta)_{\theta \in \Theta}$ a probability kernel on $\mathcal{Y}$, and $\pi$ a prior distribution on $\Theta$. Let $S$ be a sufficient statistic $S\colon \mathcal{X} \to \mathcal{Y}$ that is \emph{locally injective}, in the sense that $S(\mathbf{x}) \neq S(\mathbf{x}')$ for each pair of neighbouring inputs $\mathbf{x},\mathbf{x}' \in \mathcal{X}$ (\ie, $d(\mathbf{x},\mathbf{x}')=1$). Then, given a \emph{loss function} $L\colon \Theta \times \Theta \to \mathbb{R}_0^+$ and a {randomized private Bayes estimator} $\delta_{p,\varepsilon}$, we define the \emph{Bayes risk} of $\delta_{p,\varepsilon}$ as
\begin{multline} \label{eq:Bayes_risk}
R(\delta_{p,\varepsilon}, \pi) = \\
\int\displaylimits_{\theta \in \Theta} 
\int\displaylimits_{\mathbf{y} \in \mathcal{Y}} \int\displaylimits_{\tilde{\theta} \in \Theta} L(\theta, \tilde{\theta}) \ \delta_{p,\varepsilon}(\tilde{\theta}\vert \mathbf{y}) \, q_\theta(\mathbf{y}) \, \pi(\theta) d\theta d\mathbf{y} d\tilde{\theta}.
\end{multline}
Note that $R(\delta_{p,\varepsilon}, \pi)$ is a linear function of $\delta_{p,\varepsilon}$. In the standard Bayes point estimation, the estimate is considered as a deterministic function of the samples, but here, due to the DP constraint, the estimate will have to rely on some type of randomization mechanism. This is contrary to the standard Laplace mechanism, where the randomization is explicit. We will construct an estimator defined in terms of $\delta_{p,\varepsilon}$ that is ``optimal" in the sense of satisfying~\eqref{eq: Bayesian point estimate} subject to a DP constraint.

From the definition of DP~\eqref{eq: Differential privacy}, and that $S$ is locally injective,  it follows that for each $\mathbf{x}, \mathbf{x}' \in \mathcal{X}$ with $d(\mathbf{x},\mathbf{x}') = 1$ ($d$ is the Hamming distance), there exist $\mathbf{y}=S(\mathbf{x})$, and $\mathbf{y}'=S(\mathbf{x}')$ such that,
\begin{equation}\label{eq:diff_privacy_Bayes}
    \delta_{p,\varepsilon}(\tilde{\theta}\vert \mathbf{y}) \leq e^\varepsilon \delta_{p,\varepsilon}(\tilde{\theta}\vert \mathbf{y}'), \text{ for each } \, \tilde{\theta} \in \Theta.
\end{equation}

%

Note that \eqref{eq:diff_privacy_Bayes} is a linear constraint on $\delta_{p,\varepsilon}$.

By combining \eqref{eq:Bayes_risk} and \eqref{eq:diff_privacy_Bayes}, we can define the optimal $\varepsilon$-DP Bayes estimator $\delta_{p, \, \varepsilon}^*$ as the minimizer of the following optimization program:
%
%
%
%
%
\begin{align} \label{eq:diff_private_opt_program}
\begin{array}{cl}
\displaystyle \min_{\delta_{p,\varepsilon} \in \mathcal{P}(\mathcal{Y},\Theta)} & \displaystyle \int\displaylimits_{\Theta} \int\displaylimits_{\mathcal{Y}} \int\displaylimits_{\Theta} L(\theta, \tilde{\theta}) \delta_{p,\varepsilon}(\tilde{\theta}\vert \mathbf{y}) q_\theta(\mathbf{y}) \pi(\theta) d\theta d \mathbf{y} d\tilde{\theta} \\
\text{s.t.} & \displaystyle \delta_{p,\varepsilon}(\tilde{\theta}\vert S(\mathbf{x})) \leq e^\varepsilon \delta_{p,\varepsilon}(\tilde{\theta}\vert S(\mathbf{x}')),\text{ for each } \tilde{\theta} \in \Theta \\
&\qquad \text{and $\mathbf{x}, \mathbf{x}' \in \mathcal{X}$ s.t. $d(\mathbf{x},\mathbf{x}')=1$}, \\
 &\displaystyle \int\displaylimits_{\Theta} \delta_{p,\varepsilon}(\tilde{\theta}\vert \mathbf{y}) d\mathbf{y}  = 1, \text{for each } \mathbf{y} \in \mathcal{Y},\\
 &\displaystyle \delta_{p,\varepsilon}(\tilde{\theta}\vert \mathbf{y}) \geq 0, \text{for each } \mathbf{y} \in \mathcal{Y}, \, \tilde{\theta} \in \Theta,
\end{array} \nonumber \\
\end{align}

where $\mathcal{P}(\mathcal{Y},\Theta)$ is the set of probability densities on $\Theta$ conditioned on $\mathcal{Y}$.

\begin{figure}[t]
    \centering
    \includegraphics[width=0.8\linewidth]{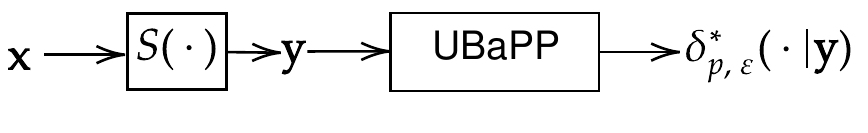}
    \caption{UBaPP estimator. Unlike the Laplace mechanism that imposes an explicit randomization, UBaPP uses an implicit randomization.}
    \label{fig: UBaPP}
\end{figure}

From \eqref{eq:diff_private_opt_program}, it is evident that $\delta_{p, \, \varepsilon}^*$  corresponds to the optimal Bayes estimator (in the sense of minimizing the Bayes risk) satisfying the DP constraint~\eqref{eq: Differential privacy}.
We denote the optimization program~\eqref{eq:diff_private_opt_program} as \emph{Unified Bayes Private Point} (UBaPP) estimator (Fig.~\ref{fig: UBaPP}) and call its solution $\delta_{p,\varepsilon}^*$ the UBaPP estimate.

In the following subsection, we specialize the UBaPP estimator for the case where the parameter space $\Theta$ and observation space $\mathcal{Y}$ are both finite, thus arriving at a tractable convex optimization program.

\subsection{Finite Case} \label{subsec: Discrete}

Consider the finite case when both $\Theta = \{ \theta_1, \dots, \theta_{|\Theta|} \}$, $\mathcal{Y} = \{ \mathbf{y}_1, \dots, \mathbf{y}_{\abs{\mathcal{Y}}} \}$ are finite. Here, $\abs{\Theta}$ and $\abs{\mathcal{Y}}$ denote the cardinality of $\Theta$ and $\mathcal{Y}$ respectively.  Then, the randomized private Bayes estimate $\delta_{p,\varepsilon}(\cdot\vert \mathbf{y})$ ($\mathbf{y}=S(\mathbf{x})$ for some $\mathbf{x} \in \mathcal{X}$) can be described in terms of a matrix $\mathbf{P} \in \mathbb{R}^{\abs{\Theta} \times \abs{\mathcal{Y}}}$ such that
\begin{equation} \label{eq: discrete_prob_kernel}
    \mathbf{P}_{i,j} := \mathbb{P}[\hat{\theta}=\theta_i \vert  \mathbf{y}_j] = \delta_{p,\varepsilon}(\theta_i \vert \mathbf{y}_j),
\end{equation}
for $i \in \{1, \dots, \abs{\Theta}\}$, $j \in \{1,\dots,\abs{\mathcal{Y}}\}$, and $\mathbf{y}_j =S(\mathbf{x}_j)$ for some $\mathbf{x}_j \in \mathcal{X}$.

Combining~\eqref{eq: discrete_prob_kernel} with~\eqref{eq:diff_private_opt_program}, we obtain the optimization program for the finite case as
\begin{align} \label{eq: discrete_case_opt_program}
     \displaystyle &\min_{\mathbf{P} \in \mathbb{R}^{\abs{\Theta} \times \abs{\mathcal{Y}}}}   \displaystyle \sum_{i=1}^{\abs{\mathcal{Y}}} \sum_{j=1}^{\abs{\Theta}} \sum_{k=1}^{\abs{\Theta}} L(\theta_j, \theta_k) \mathbf{P}_{k,i} \mathbb{P}[\mathbf{y}_i \vert \theta = \theta_j] \pi(\theta_j) \nonumber \\
   &\quad \quad \text{s.t.}  \quad \,\, \displaystyle \mathbf{P}_{k,i} \leq e^\varepsilon \mathbf{P}_{k,i'}, \,  \text{for all } k \in \{1, \ldots, \abs{\Theta}\}  \nonumber \\ 
   &\quad \quad \qquad \qquad \text{and } i,i' \in \{1,\ldots,\abs{\mathcal{Y}}\} \text{ s.t. } d(\mathbf{x}_i,\mathbf{x}_{i'}) = 1, \nonumber\\
   &\quad \quad \quad \quad \,\,\,\, \mathbf{1}^T\mathbf{P} = \mathbf{1}^T, \nonumber\\
   &\quad \quad \quad \quad \,\,\,\, \mathbf{P} \geq 0,
\end{align}
where $\pi\colon \Theta \to [0,1]$ is a prior probability mass function (pmf) on $\Theta$, and $\mathbf{1} = (1,\ldots,1)^T \in \mathbb{R}^{\abs{\Theta}}$. Here, $\mathbf{P} \geq 0$ means that all the entries of $\mathbf{P}$ are non-negative. Note that $\mathbb{P}[\mathbf{y}_i \vert \theta = \theta_j] = q_{\theta_j}(\mathbf{y}_i)$.

By defining $\mathbf{L} \in \mathbb{R}^{\abs{\Theta} \times \abs{\Theta}}$ 
as $\mathbf{L}_{j,k} := L(\theta_j, \theta_k)$
for $j,k \in \{1,\ldots,\abs{\Theta}\}$, $\mathbf{Q} \in \mathbb{R}^{\abs{\mathcal{Y}} \times \abs{\Theta}}$ as $\mathbf{Q}_{i,j} := \mathbb{P}[\mathbf{y}_i | \theta = \theta_j]=q_{\theta_j}(\mathbf{y}_i)$, for  $i \in \{1,\ldots,\abs{\mathcal{Y}}\}$, $j \in \{1,\ldots,\abs{\Theta}\}$, and $\boldsymbol{\pi} \in \mathbb{R}^{\abs{\Theta}}$ as $\boldsymbol{\pi}_j := \pi(\theta_j)$ ($j \in \{1,\ldots,\abs{\Theta}\}$), we can re-write the optimization problem as
\begin{align} \label{eq: final_opt_program_discrete}
\displaystyle &\min_{\mathbf{P} \in \mathbb{R}^{\abs{\Theta} \times \abs{\mathcal{Y}}}}  \displaystyle \text{tr}(\mathbf{Q} \text{diag}(\boldsymbol{\pi}) \mathbf{L} \mathbf{P}) \nonumber \\
&\quad \quad \text{s.t.}  \quad \,\, \displaystyle \mathbf{P}_{k,i} \leq e^\varepsilon \, \mathbf{P}_{k,i'}, \,  \text{for all } k \in \{1, \ldots, \abs{\Theta}\}, \nonumber \\ 
   & \quad \quad \quad \quad \quad \,\,\,\, \text{ and } i,i' \in \{1,\ldots,\abs{\mathcal{Y}}\} \text{ s.t. } d(\mathbf{x}_i,\mathbf{x}_{i'}) = 1 \nonumber,\\
   &\quad \quad \quad \quad \,\,\,\, \mathbf{1}^T\mathbf{P} = \mathbf{1}^T, \nonumber\\
   &\quad \quad \quad \quad \,\,\,\, \mathbf{P} \geq 0.
\end{align}
Here, $\text{tr}(\cdot)$ denotes the trace of a matrix and $\text{diag}(\boldsymbol{\pi})$ denotes a diagonal matrix whose $j^{th}$ entry is $\boldsymbol{\pi}_j = \pi(\theta_j)$ ($j \in \{1,\ldots,\abs{\Theta}\}$). The minimizer of optimization program~\eqref{eq: final_opt_program_discrete} gives us the UBaPP estimate $\delta_{p, \varepsilon}^{*}$ for the finite case, and it can be obtained using CVXPY~\citep{diamond2016cvxpy}. To evaluate the performance of the UBaPP estimator, we compute the theoretical mean-square error (MSE) of $\delta_{p, \varepsilon}^{*}$ for different values of $\varepsilon$.  We summarize the computation of $\delta_{p, \varepsilon}^{*}$ and its MSE in Algorithm~\ref{alg: UBaPP}.

%
\begin{algorithm}
  \caption{UBaPP estimator}

\begin{algorithmic}[1]
 \State \textbf{Input}: $\mathcal{X}$, $\varepsilon > 0$, $\boldsymbol{\pi}$, and $\Theta$
 \State Compute $\mathbf{L}$ by evaluating $L(\theta_j,\theta_k)$ for $j,k \in \{1,\dots,\abs{\Theta}\}$
 \State Compute $\mathbf{Q}$ by evaluating $\mathbb{P}[\mathbf{y}_i \vert \theta = \theta_j]$ for each $i \in \{1,\ldots, \abs{\mathcal{Y}}\}$, $j \in \{1,\ldots,\abs{\Theta}\}$
 \State Solve \eqref{eq: final_opt_program_discrete} by CVXPY solver and store the minimizer as $\delta_{p,\varepsilon}^*$
  
 \State \textbf{Output}: $\delta_{p,\varepsilon}^*$, $\text{MSE}=\text{tr}(\mathbf{Q} \text{diag}(\boldsymbol{\pi}) \mathbf{L} \delta_{p,\varepsilon}^*)$. 
 \end{algorithmic} 
 \label{alg: UBaPP}
\end{algorithm}




Finally, we compare the MSE of the UBaPP estimator with that of the Laplace Private Bayes estimator $\delta_{lpb,\varepsilon}$ \eqref{eq: Laplace-private-Bayes}, based on the Laplace mechanism, whose MSE is computed as in Algorithm~\ref{alg: Laplace_MSE}.

\begin{algorithm}
  \caption{MSE of LBaPP estimator}

\begin{algorithmic}[1]
 \State \textbf{Input}: $\mathcal{X}$, $\varepsilon > 0$, $\boldsymbol{\pi}$, $N$ (number of Monte-Carlo runs), and $\Theta$
 \State \textit{Initialize}: cum\_mse = $0$
 \For{$i=1,\ldots,N$}
 \State Sample $\theta_i \sim \boldsymbol{\pi}$
 \State Sample $\mathbf{x}_i \sim \mathbb{P}^n(\cdot\vert\theta_i)$
 \State Compute $\delta^*_\text{Bayes}(\mathbf{x}_i) $ from \eqref{eq: Bayesian point estimate}
 \State Compute $\sigma_{\delta_{\text{Bayes}}^*}$ from \eqref{eq: l1_sensitivity_Bayes} 
 \State $\hat{\theta}_i = \pb(\mathbf{x}_i) = \delta_{\text{Bayes}}^*(\mathbf{x}_i) + \text{Lap}\left(0,\dfrac{\sigma_{\delta_{\text{Bayes}}^*}}{\varepsilon} \right)$
 \State cum\_mse $\leftarrow$ cum\_mse + $(\hat{\theta}_i -\theta_i)^2$
 \EndFor
 \State \textbf{Output}: MSE $=\dfrac{\text{cum\_mse}}{N}$.
\end{algorithmic} 
 \label{alg: Laplace_MSE}
\end{algorithm}

\section{Simulations} \label{sec: Simulation}
In this section, we compare our approach with the standard Laplace mechanism. In particular, we demonstrate that the accuracy-privacy trade off for our approach is better than for the Laplace mechanism, especially in the high privacy regime (\ie, for small values of $\varepsilon$). For this purpose, we consider the simple numerical example outlined next.

\subsection{Setup}\label{subsec: setup_sim}
Suppose we are interested in estimating the unknown parameter $\theta \in \Theta = [0,1]$ of a Bernoulli distribution $\text{Ber}(\theta)$ from which $K$ independent samples are generated in $\{0,1\}$. The simulation setup is as follows:
\begin{itemize}
\item Let $\pi$ be the uniform distribution over $\Theta$.
\item Let $x_1,\ldots,x_K \stackrel{i.i.d.}{\sim} \text{Ber}(\theta_i)$ be outcomes of Bernoulli trials, where $\theta_i \in \Theta$, and $\mathbf{x}=(x_1,\ldots,x_K)^T$. Then, define the sufficient statistic $S$ as
\begin{equation*}
    \mathbf{y}=S(\mathbf{x}) = \sum\limits_{i=1}^{K} x_i.
\end{equation*}
This means that $\mathbf{y} \in \mathcal{Y}=\{0,\ldots,K\}$.
\item Let $L$ be the square loss, \ie, $L(\theta_i,\theta_j) = (\theta_i-\theta_j)^2$. 
\end{itemize}

For this setup, we first derive the LBaPP estimator $\delta_{lpb,\varepsilon}$ in Subsection~\ref{subsec: sim_lpb}. Then, in Subsection~\ref{subsec: sim_rpb}, we describe the procedure to obtain the UBaPP estimator $\delta_{p,\varepsilon}^*$.

\subsection{LBaPP estimator}\label{subsec: sim_lpb}
Let $\mathbb{P}(\theta \vert \mathbf{y})$ denote the posterior distribution of $\theta$. Then, it can be shown that
\begin{equation}\label{eq: beta_post}
\mathbb{P}(\theta\vert \mathbf{y}) \propto \theta^{\left(\sum_{i=1}^{K} x_i +1 \right) -1} (1-\theta)^{\left(K+1-\sum_{i=1}^{K}x_i\right) -1}.
\end{equation}
This implies that $\mathbb{P}(\theta\vert \mathbf{y})$ is a Beta distribution. Therefore, from \eqref{eq: conditional_mean_estimate}, the non-private Bayes estimate is given by
\begin{equation}\label{eq: non_private_Bayes_sim}
    \delta^*_{\text{Bayes}}(\mathbf{x}) = \dfrac{\sum\limits_{i=1}^{K} x_i+1}{K+2} .
\end{equation}
We now add Laplace noise to obtain $\delta_{lpb,\varepsilon}(\mathbf{x})$. To this end, we need to compute $\sigma_{\delta_{\text{Bayes}}^*}$ using~\eqref{eq: l1_sensitivity_Bayes}. Let $\mathbf{x}=(x_1,\ldots,x_K)^T$, and $\mathbf{x}'=(x'_1,\ldots,x'_K)^T$. Then,
\begin{align}\label{eq: l1_sensitivity_Bayes_sim}
    \sigma_{\delta_{\text{Bayes}}^*} &= \sup{\mathbf{x},\mathbf{x}': \, d(\mathbf{x},\mathbf{x}')=1} \, \, \vert \vert \delta_{\text{Bayes}}^*(\mathbf{x}) - \delta_{\text{Bayes}}^*(\mathbf{x}') \vert \vert_1 \nonumber \\
    &\stackrel{}{=} \dfrac{1}{K+2} \, \, \sup{\mathbf{x},\mathbf{x}': \, d(\mathbf{x},\mathbf{x}')=1} \abs{\sum\limits_{i=1}^K x_i - \sum\limits_{i=1}^K x'_i} \nonumber \\
&= \dfrac{1}{K+2}.
\end{align}
Finally, using~\eqref{eq: l1_sensitivity_Bayes_sim} and \eqref{eq: non_private_Bayes_sim} in~\eqref{eq: Laplace-private-Bayes}, we obtain 
\begin{equation} \label{eq: LBaPP_Bernoulli}
    \delta_{lpb,\varepsilon}(\mathbf{x}) = \dfrac{\sum_{i=1}^{K} x_i+1}{K+2} + \text{Lap}\left(0,\dfrac{1}{(K+2) \varepsilon} \right).
\end{equation}

We compute the MSE of LBaPP estimator using Algorithm~\ref{alg: Laplace_MSE}, where we take $N=5000$.

\subsection{UBaPP estimator}\label{subsec: sim_rpb}
To obtain the UBaPP estimator $\delta^*_{p,\varepsilon}$ we will use formulation~\eqref{eq: final_opt_program_discrete}, for which we only need to discretize $\Theta$, as $\mathcal{Y}$ is already finite.
For this purpose, we consider a grid of $M_{\theta}$ equally spaced points and denote the set of such equally spaced points by $\tilde{\Theta}$.
Due to the discretization, we note that there should be an additional factor $1 / M_{\theta}$ in the objective function of \eqref{eq: final_opt_program_discrete}, but this factor can be ignored since it does not depend on $\delta_{p,\varepsilon}$. Also, $\pi$ is now the uniform distribution over $\tilde{\Theta}$, and hence, the $j^{th}$ entry  of $\text{diag}(\boldsymbol{\pi})$ is $\boldsymbol{\pi}_j = 1/M_{\theta}$ ($j \in \{1,\ldots,|\tilde{\Theta}|\}$). For this simulation setup, it is easy to see that $\mathbf{Q}$ is Binomial, \ie,
\begin{equation*}
    \mathbf{Q}_{i,j} = \mathbb{P}[\mathbf{y}_i| \theta_j] = {K \choose \mathbf{y}_i} \theta_j^{\mathbf{y}_i} (1-\theta_j)^{K-\mathbf{y}_i},
\end{equation*}
where $\theta_j \in \tilde{\Theta}$, and $\mathbf{y}_i \in \mathcal{Y} = \{0,\ldots,K\}$, for $j \in \{1,\ldots,|{\tilde{\Theta}}|\}$, and $i \in \{0,\ldots,\abs{\mathcal{Y}}\}$ respectively. Also, for the simulation, we take $M_{\theta} =5000$. 

\begin{figure}[httb!]
\centering
\includegraphics[width=1\linewidth]{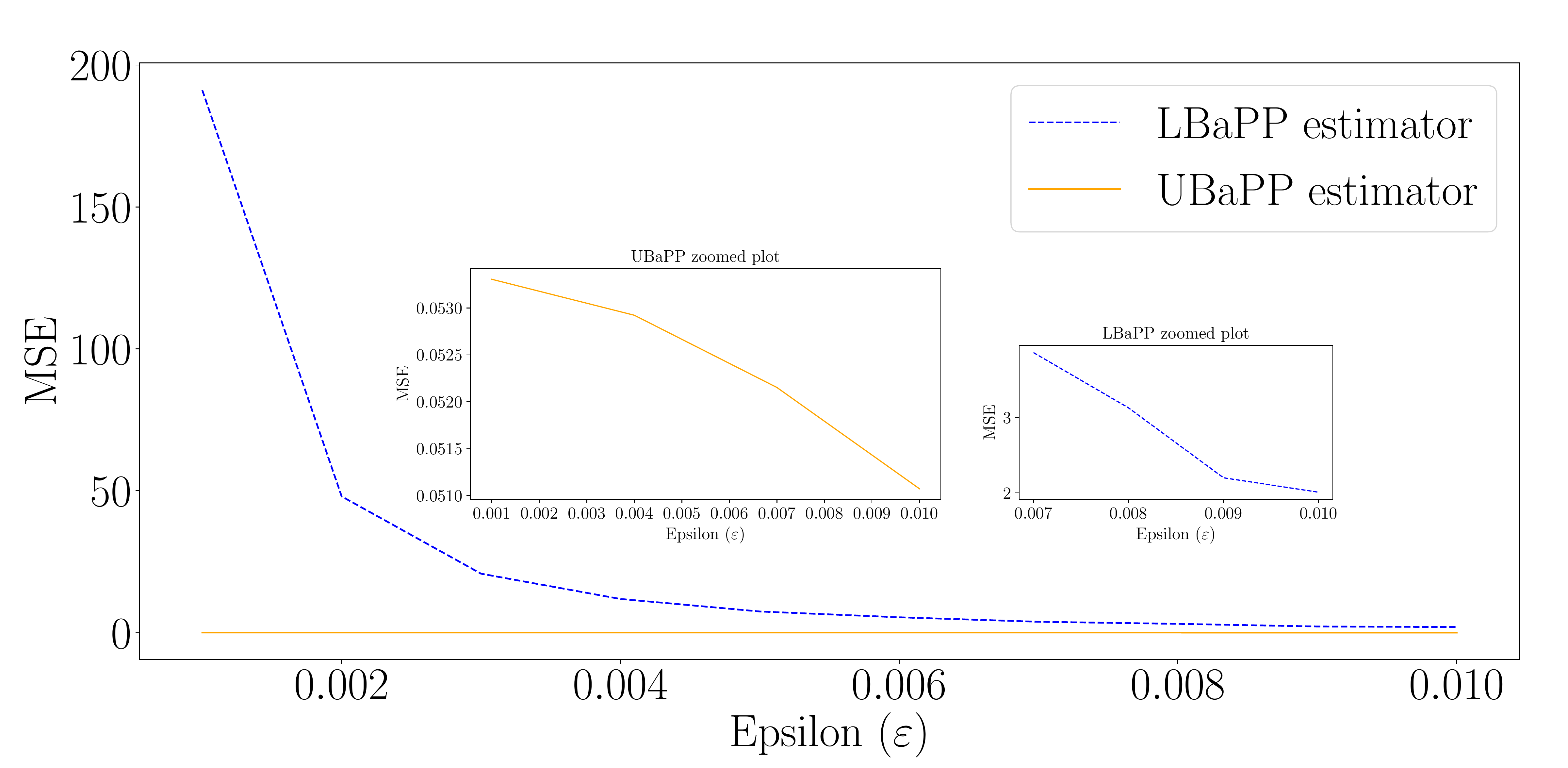}
\caption{Plot of MSE vs. $\varepsilon$ for the LBaPP and UBaPP estimators in the high privacy regime. A zoomed plot for the UBaPP and LBaPP estimators is included.}
%
\label{fig: high_privacy_regime}
\end{figure}
    
\begin{figure}[httb!]
\centering
\includegraphics[width=1\linewidth]{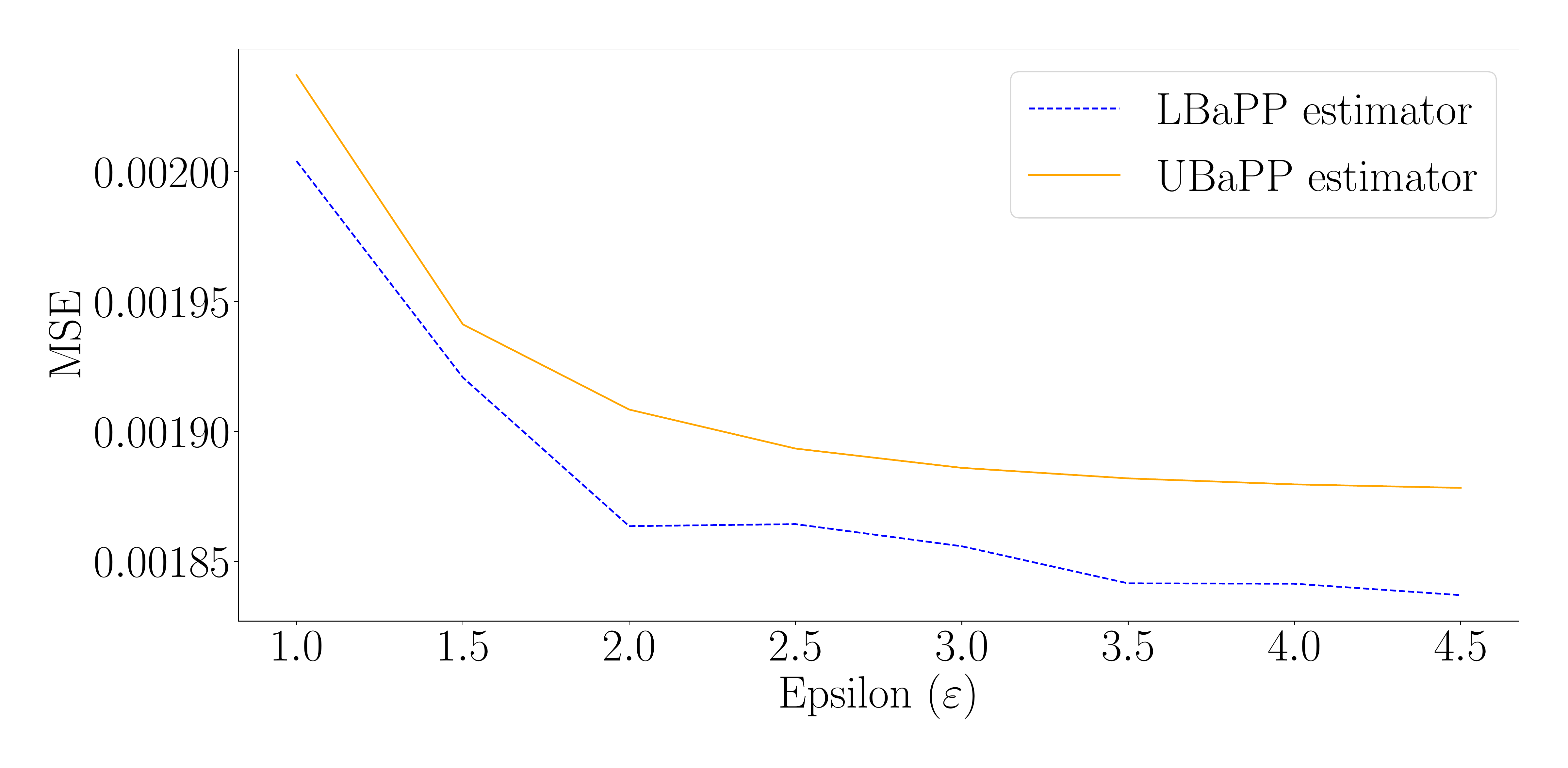}
\caption{Plot of MSE vs. $\varepsilon$ for the LBaPP and UBaPP estimators in the moderate-low privacy regime.}
\label{fig: low_privacy_regime}
\end{figure}

\subsection{Plots}
We now study the performance of using our approach to differentially private Bayes point estimation, UBaPP.
%
In particular, we focus on 
\begin{itemize}
    \item the effect of the privacy parameter $\varepsilon$ on the MSE, for a fixed number of Bernoulli trials;
    \item the effect of the number of Bernoulli trials ($K$) on the MSE, for a fixed $\varepsilon$.
\end{itemize}

To analyze the MSE of these two estimators for different values of $\varepsilon$, we consider two ranges of $\varepsilon$, which are classified as ``High privacy regime'' and ``Moderate-low privacy regime''. The high privacy regimes corresponds to low values of $\varepsilon$, \ie, to $\varepsilon  \in \{0.001,0.002,\ldots,0.010\}$, while the moderate-low privacy regimes corresponds to high values of $\varepsilon$, \ie, $\varepsilon \in \{1,1.5,\ldots,5\}$. The resulting plots of the MSE vs. $\varepsilon$ of the LBaPP and UBaPP estimators for the high and moderate-low privacy regimes are shown in Figs.~\ref{fig: high_privacy_regime} and~\ref{fig: low_privacy_regime}, respectively.



\begin{figure}[httb!]

\includegraphics[width=1\linewidth]{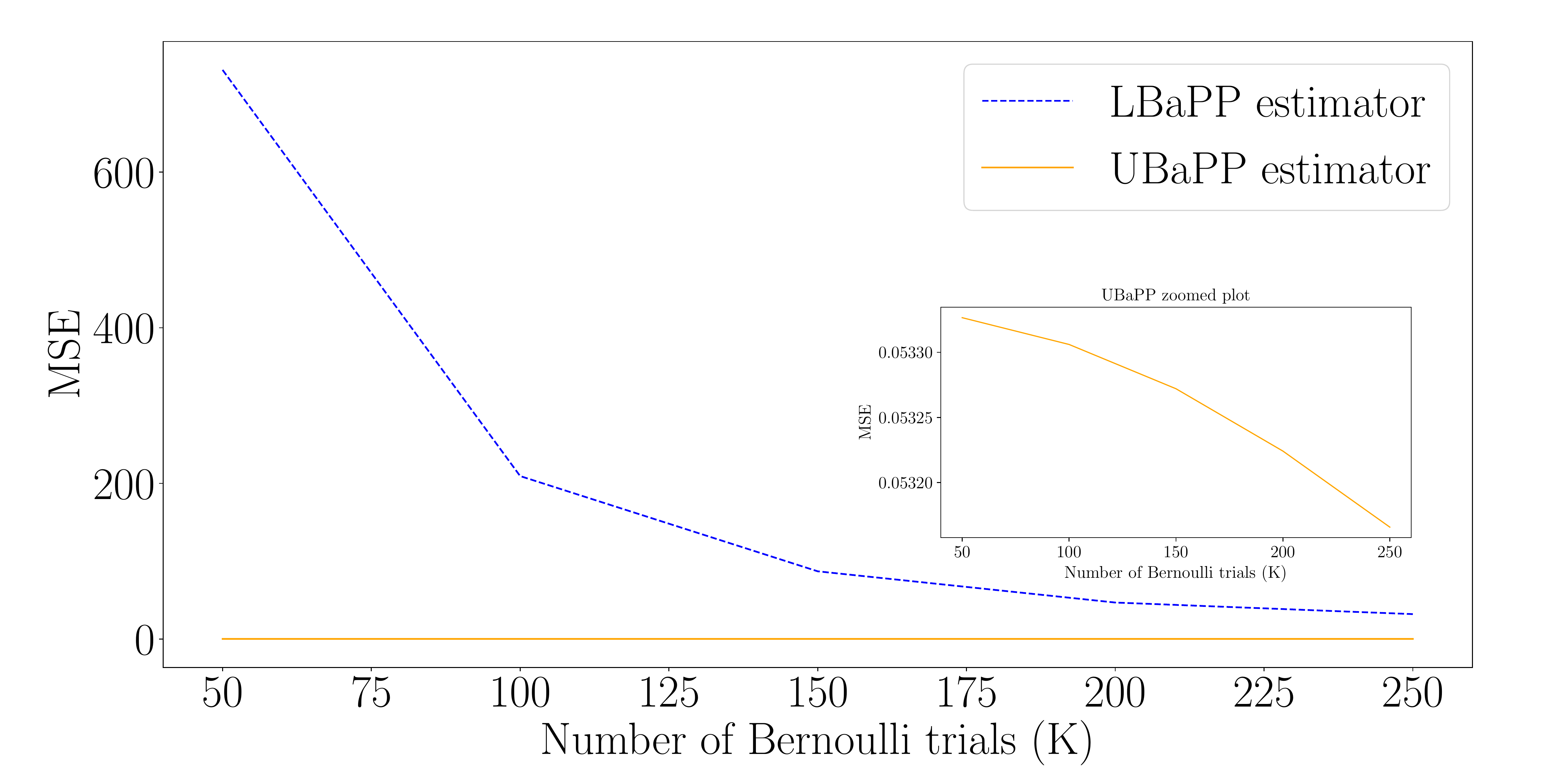}
\caption{Plot of MSE vs. number of Bernoulli trials (K) for the LBaPP and UBaPP estimators in a fixed high privacy regime. A zoomed plot for the UBaPP estimator is included.}
\label{fig: high_privacy_regime_diff_trials}
\end{figure}

\begin{figure}[httb!]
\includegraphics[width=1\linewidth]{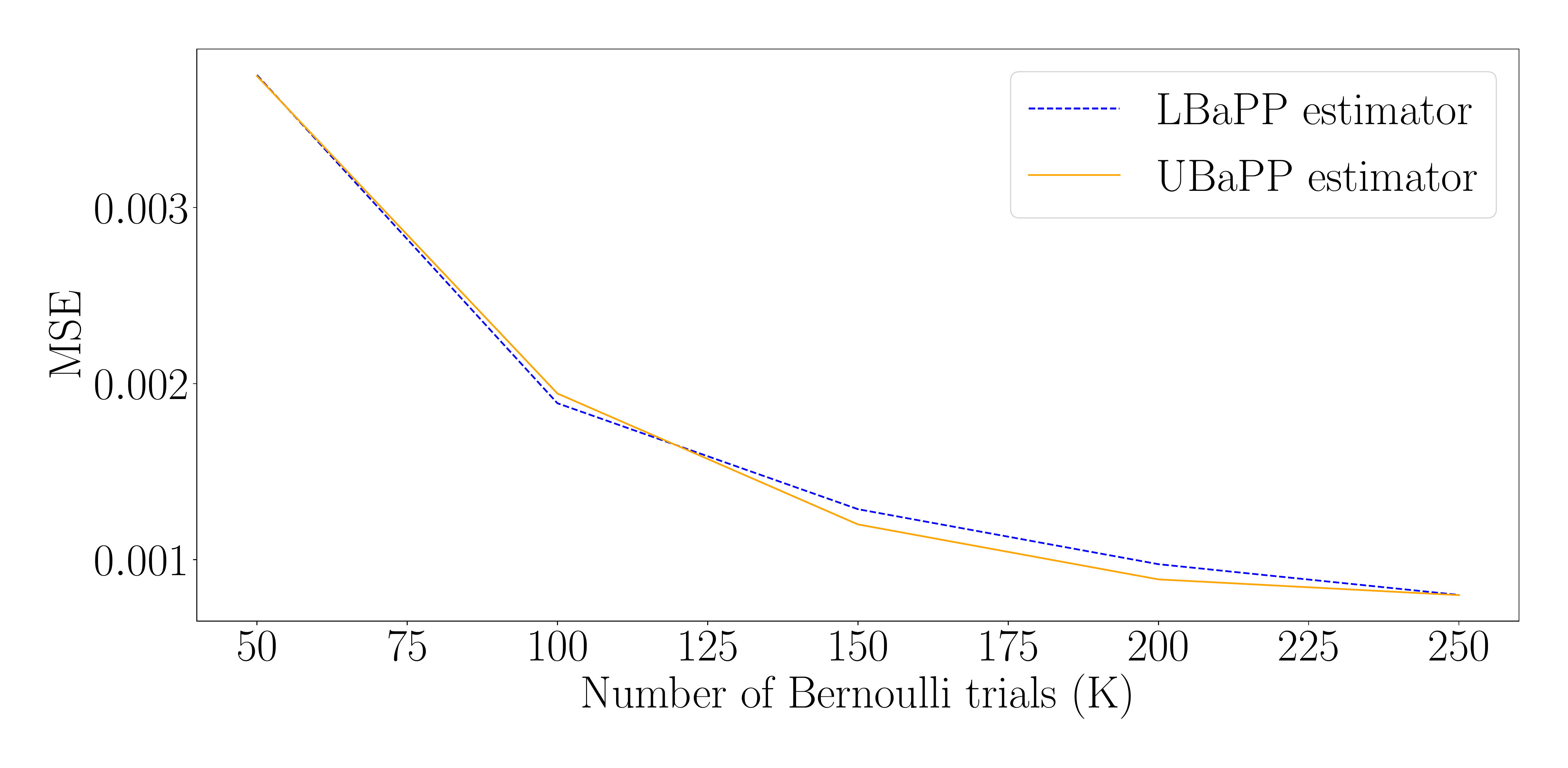}
\caption{Plot of MSE vs. number of Bernoulli trials (K) for the LBaPP and UBaPP estimators in a fixed moderate-low privacy regime.}
\label{fig: low_privacy_regime_diff_trials}
\end{figure}

Next, to study the effect of the number of Bernoulli trials on the MSE of both estimators, we plot MSE of LBaPP and UBaPP for different values of $K$. We consider two different values of $\varepsilon$, namely, $\varepsilon=0.001$ and $\varepsilon=5$, corresponding to the high privacy and moderate-low privacy regimes, respectively. The corresponding plots are shown in Figs.~\ref{fig: high_privacy_regime_diff_trials} and~\ref{fig: low_privacy_regime_diff_trials}, respectively.

\begin{figure*}[httb!]
\begin{subfigure}{0.5\textwidth}
\includegraphics[width=1\linewidth]{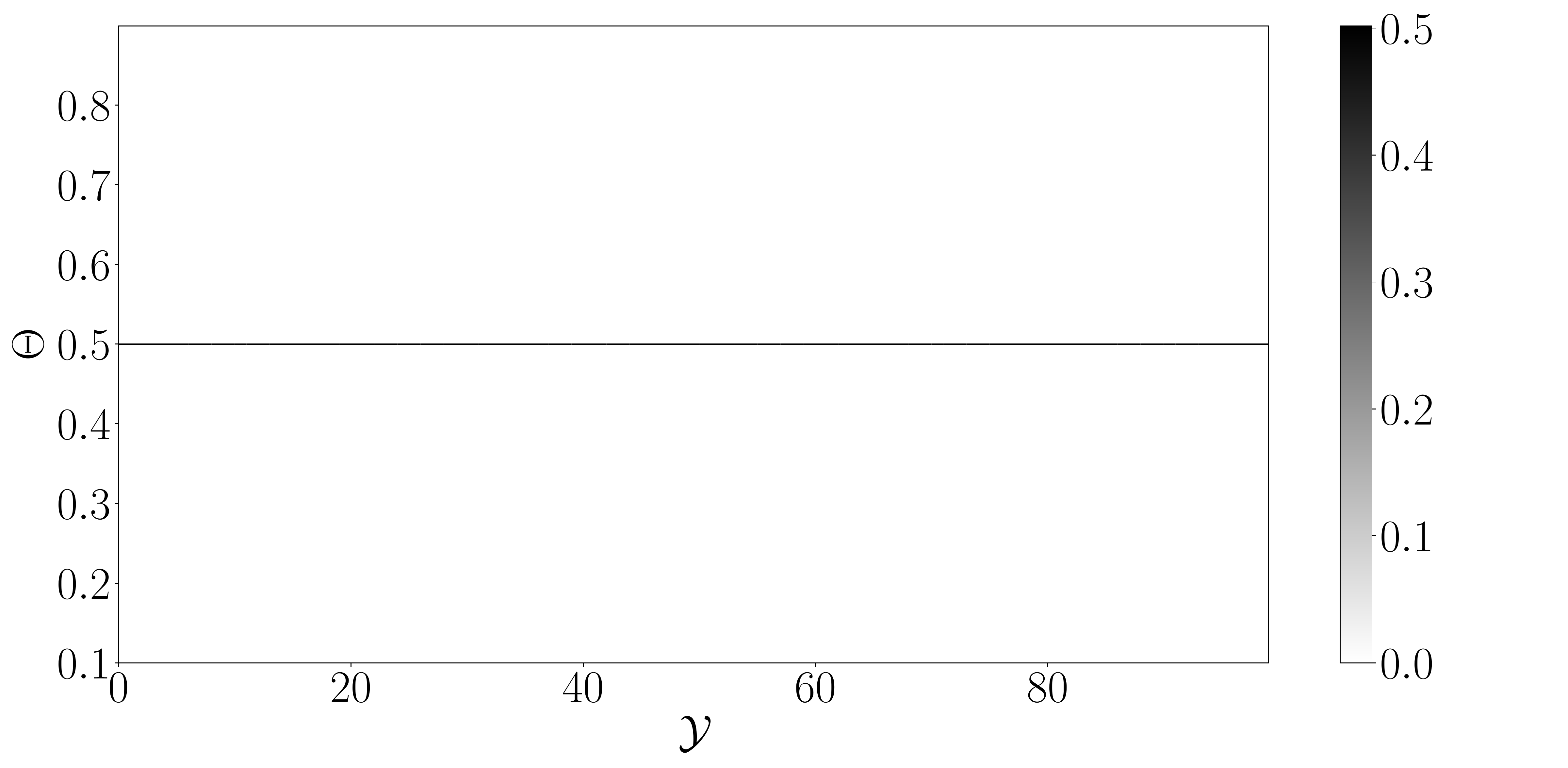}
\caption{$\varepsilon=0.0001$}
\label{fig: heat_map_eps_0.0001}
\end{subfigure}
\hfill
\begin{subfigure}{0.5\textwidth}
\includegraphics[width=1\linewidth]{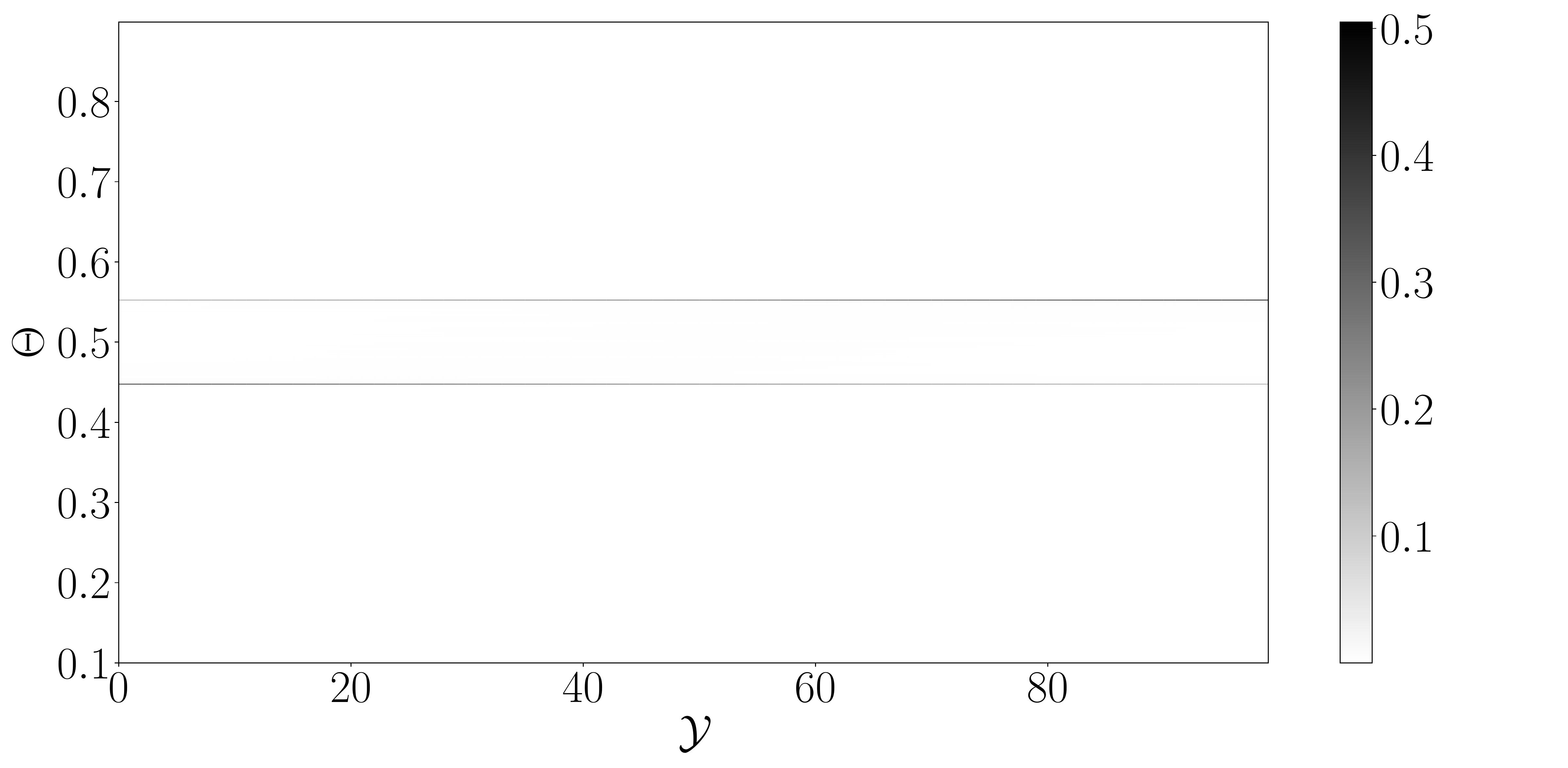}
\caption{$\varepsilon=0.01$}
\label{fig: heat_map_eps_0.01}
\end{subfigure}
\vfill
\begin{subfigure}{0.5\textwidth}
\includegraphics[width=1\linewidth]{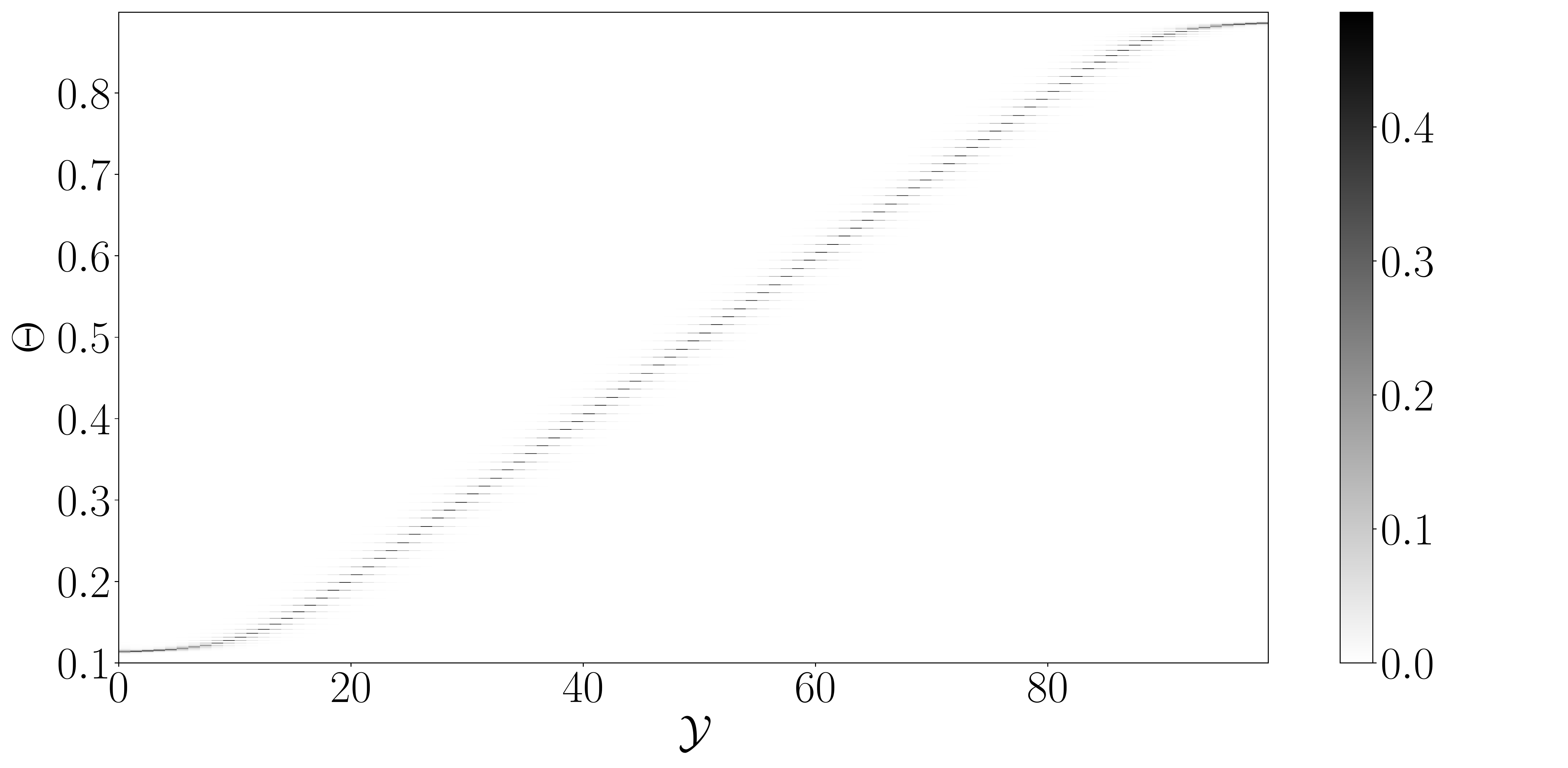}
\caption{$\varepsilon=1$}
\label{fig: heat_map_eps_1}
\end{subfigure}
\hfill
\begin{subfigure}{0.5\textwidth}
\includegraphics[width=1\linewidth]{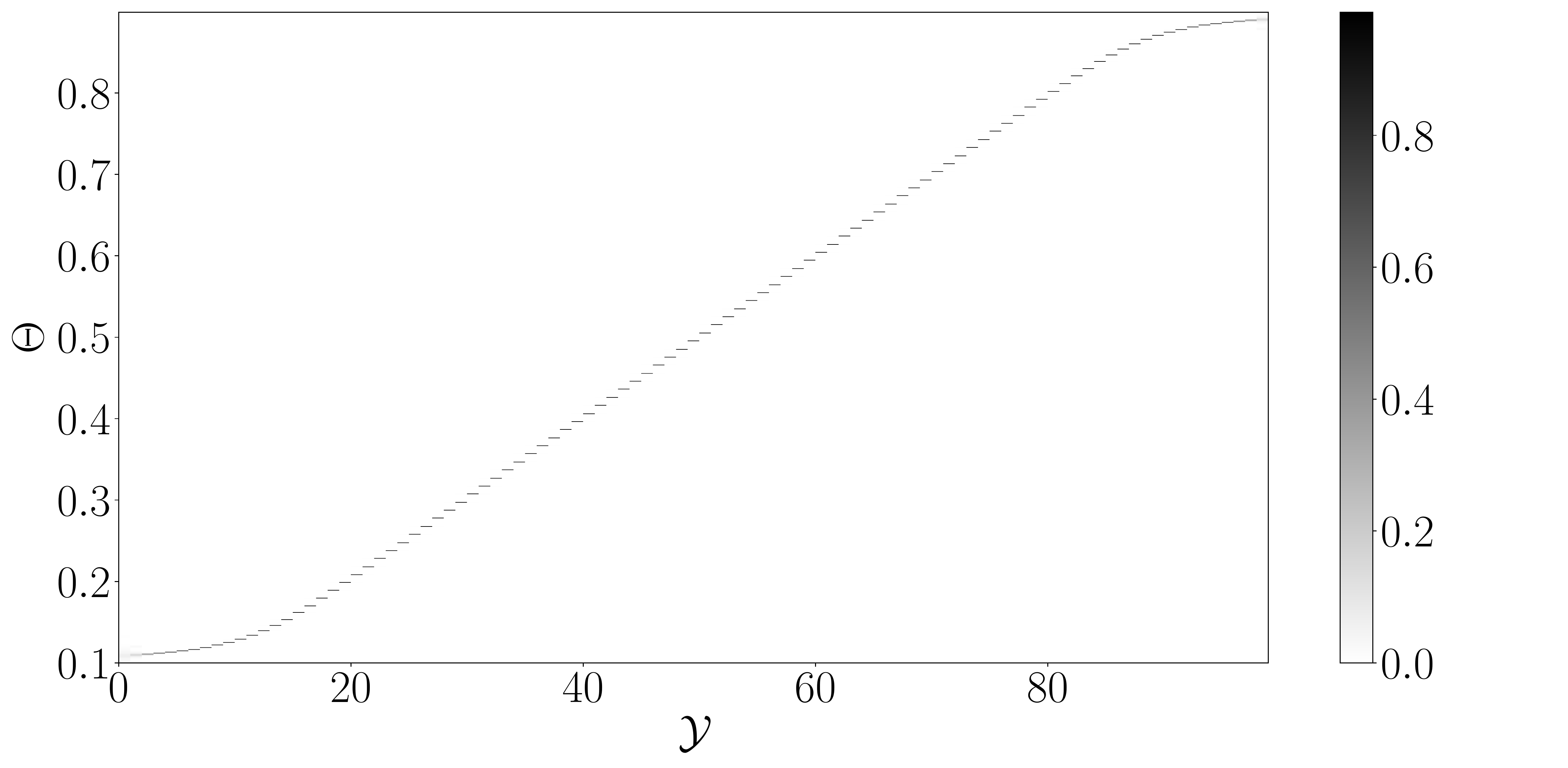}
\caption{$\varepsilon=5$}
\label{fig: heat_map_eps_5}
\end{subfigure}
\caption{\centering Heat maps of UBaPP estimate $\delta_{p,\varepsilon}^*$ for different values of $\varepsilon$.}
\label{fig: heat_maps}
\end{figure*}

Finally, to understand the randomization induced by the UBaPP estimator, we plot in Fig.~\ref{fig: heat_maps} heat maps of UBaPP estimate $\delta_{p,\varepsilon}^*$ for different values of $\varepsilon$.

\subsection{Discussion}
We observe from Fig.~\ref{fig: high_privacy_regime} that UBaPP provides high accuracy (\ie, low MSE) in the high privacy regime, when compared to that of LBaPP.
Also, Fig.~\ref{fig: low_privacy_regime} shows that in the moderate-low privacy regime, the accuracy of UBaPP is very similar to that of LBaPP estimator: even though the MSE of UBaPP is higher than that of LBaPP, the difference is only of order approx. $10^{-3}$. This difference is due to the discretization of the parameter space $\Theta$ used in the implementation of UBaPP, whereas the computation of the LBaPP estimates is exact (\ie, no discretization is employed).

Regarding the dependence on the number of trials $K$, from Fig.~\ref{fig: high_privacy_regime_diff_trials} it is clear that in the high privacy regime, the MSE of UBaPP is less dependent on $K$ than the MSE of LBaPP. This means that, for low values of $K$, UBaPP has significantly higher accuracy than LBaPP. On the other hand, according to Fig~\ref{fig: low_privacy_regime_diff_trials}, the accuracy of LBaPP and UBaPP is similar in the moderate-low privacy regime.

To understand the difference between our approach, UBaPP, and the Laplace mechanism LBaPP, notice from Fig.~\ref{fig: heat_map_eps_0.0001} that, for low values of $\varepsilon$, UBaPP outputs a deterministic estimate around $\theta = 0.5$ (which corresponds to the mean of the prior distribution on $\theta$) for every ${y} \in \mathcal{Y}$; since low values of $\varepsilon$ imply high privacy, no valid inference about the data (input) $\mathbf{x}$ can be made based on the estimate of $\theta$, so the estimator becomes independent of the data by outputting the same deterministic estimate for all ${y} \in \mathcal{Y}$.
Also, as seen in Figs.~\ref{fig: heat_map_eps_0.01}-\ref{fig: heat_map_eps_1}, when $\varepsilon$ is increased (which implies a shift from high to moderate levels of privacy), UBaPP introduces some level of randomization, because for a given value of $y \in \mathcal{Y}$ the probability distribution  $\delta_{p,\varepsilon}^*$ is not concentrated at a single value of $\theta \in \Theta$. Finally, as $\varepsilon$ is further increased to very low privacy levels, Fig.~\ref{fig: heat_map_eps_5} shows that UBaPP becomes deterministic again, as it tends to the standard (non-private) Bayes point estimator of $\theta$, which is known to be deterministic for convex loss functions~\citep{Berger-85}. Also, note from Fig.~\ref{fig: heat_map_eps_5} that the deterministic estimate varies with $y \in \mathcal{Y}$, which suggests that UBaPP estimator becomes strongly dependent of the data, and hence it becomes non-private.

In contrast, the Laplace mechanism adds randomization to the standard non-private Bayes point estimator of $\theta$, with a variance that increases as $\varepsilon \to 0$. Thus, for large values of $\varepsilon$ it coincides, like UBaPP, with the standard (non-private) Bayes point estimator of $\theta$ (shown in Fig.~\ref{fig: heat_map_eps_5}). However, for very small values of $\varepsilon$ (\ie, very high privacy), the variance of the estimator grows unbounded, which implies that its MSE tends to infinity, as shown in Fig.~\ref{fig: high_privacy_regime}, whereas the MSE of UBaPP tends to a constant as $\varepsilon \to 0$.

In conclusion, we see from the numerical study that UBaPP is more accurate than LBaPP in the high privacy regime, while for moderate-low privacy constraints both estimators yield similar performance (save for the discretization needed to implement UBaPP). 


\section{Conclusion} \label{sec: Conclusion}
In this paper, we have studied the problem of Bayes point estimation under differential privacy.
We have argued that the standard approach based on the Laplace mechanism may not give accurate estimates under high privacy constraints.
We then proposed an optimal approach that combines risk minimization (minimum MSE) and differential privacy into a single convex optimization program, and specialized this approach to the case of finite parameter and observation space.
Via a simple numerical study, we have shown that our approach yields more accurate estimates in the high privacy regime than the Laplace mechanism, and that both approaches have similar performance under low privacy constraints.

As future work, we plan the extension of our approach to continuous (and even potentially high dimensional) parameter and observation spaces.



\bibliography{ifacconf}             
                                                   







\end{document}